\documentclass[11pt, reqno]{amsart}
\usepackage{amsmath}
\usepackage{amssymb}
\input amssym.def
\input amssym.tex

\font\tengothic=eufm10
\font\sevengothic=eufm7
\newfam\gothicfam
      \textfont\gothicfam=\tengothic
      \scriptfont\gothicfam=\sevengothic
\def\goth#1{{\fam\gothicfam #1}}

\numberwithin{equation}{section}

\begin{document}
\setlength{\baselineskip}{1.5em} \setcounter{section}{-1}

{\theoremstyle{plain}
	\newtheorem{thm}{\bf Theorem}[section]
	\newtheorem{pro}[thm]{\bf Proposition}
	\newtheorem{claim}[thm]{\bf Claim}
	\newtheorem{lemma}[thm]{\bf Lemma}
	\newtheorem{cor}{\bf Corollary}[thm]
}
{\theoremstyle{remark}
	\newtheorem{remark}[thm]{\bf Remark}
	\newtheorem{example}[thm]{\bf Example}
}
{\theoremstyle{definition}
	\newtheorem{mydef}[thm]{\bf Definition}
}

\newcommand{\pp}{{\mathbb P}}
\newcommand{\tx}{\tilde{X}}
\newcommand{\ti}{{\mathcal I}}
\newcommand{\dte}{{\mathcal D}_{e,t}}
\newcommand{\vte}{\tx_{e,t}}
\newcommand{\proj}{\mbox{\rm Proj }}
\newcommand{\reg}{\mbox{\rm reg }}
\newcommand{\sheaf}[1]{\mathcal #1}
\newcommand{\isomap}{\stackrel{\sim}{\rightarrow}}
\newcommand{\zz}{{\mathbb Z}}
\newcommand{\nn}{{\mathbb N}}
\newcommand{\x}{D}
\newcommand{\h}{{\mathbb H}}
\newcommand{\cc}{C}
\newcommand{\ix}{I_{\mathbb X}}
\newcommand{\cv}{{\mathcal V}}
\newcommand{\calo}{{\mathcal O}}
\newcommand{\cl}{{\mathcal L}}
\newcommand{\cf}{{\mathcal F}}
\newcommand{\ck}{{\mathcal K}}
\newcommand{\ci}{{\mathcal I}}
\newcommand{\cm}{{\mathcal M}}
\newcommand{\R}{{\mathcal R}}
\newcommand{\sfrac}[2]{\frac{\displaystyle #1}{\displaystyle #2}}

\title{Projective embeddings of projective schemes blown up at subschemes}
\author{Huy T\`ai H\`a}
\address{Department of Mathematics, University of Missouri, Columbia MO 65211}
\email{tai@math.missouri.edu}

\begin{abstract}
Suppose $X$ is a nonsingular projective scheme, $Z$ a
nonsingular closed subscheme of $X$. Let $\tx$ be the blowup of $X$ centered at $Z$,
$E_0$ the pull-back of a general hyperplane in $X$, and $E$ the exceptional
divisor. In this paper, we study projective embeddings of $\tx$ given by divisors $\dte = tE_0 - eE$. When $X$ satisfies a necessary condition, we give explicit values of $d$ and $\delta$ such that
for all $e > 0$ and $t > ed + \delta$, $\dte$ embeds $\tx$ as a projectively normal
and arithmetically Cohen-Macaulay scheme. We also give a uniform bound for the
regularities of the ideal sheaves of these embeddings, and study their asymptotic
behaviour as $t$ gets large compared to $e$. When $X$ is a surface and $Z$ is a
0-dimensional subscheme, we further show that these embeddings possess property
$N_p$ for all $t \gg e > 0$.
\end{abstract}
\maketitle

\begin{center}
{\it Dedicated to the sixtieth birthday of Prof. A.V. Geramita}
\end{center}

\section{Introduction}

Let $R$ be a finitely generated standard graded $k$-algebra of dimension $(n+1)$, and $X = \mbox{Proj } R$ a nonsingular projective scheme ($\dim X = n$). Suppose $Z$ is a nonsingular closed subscheme of $X$, and $I \subseteq R$ its defining ideal ($I$ saturated). Let $\pi: \tx \rightarrow X$ be the blowup of $X$ centered at $Z$. For positive integers $e, t \in \nn$, consider the divisor $\dte = tE_0 - eE$ on $\tx$, where $E_0$ is the pull-back to $\tx$ of a general hyperplane in $X$ and $E$ is the exceptional
divisor of the blowing up. $\dte$ corresponds to the linear system of hypersurfaces
of degree $t$  containing $Z$ with multiplicity at least $e$. For $t
\gg e > 0$, it is known (cf. \cite[Lemma 1.1]{ch}) that $\dte$ is very ample on
$\tx$. This gives a projective embedding $\tx \hookrightarrow
\pp^{N_{e,t}}$. Let $\vte$ be the image of $\tx$ in this embedding (i.e. $\tx
\isomap \vte \subseteq \pp^{N_{e,t}}$). In the last fifteen years, there  have been
many studies on these projective embeddings of $\tx$ in various situations
depending upon the scheme $X$, its subscheme $Z$, and the values of $e$ and $t$. For
instance, \cite{ch, ge-gi, ge-gi-h, ge-gi-p, gi1, gi-lo, tai-thesis, tai, holay1, holay2}. This line of works also has a close relation with the studies of diagonal subalgebras of a bi-graded algebra, such as in \cite{chtv, stv, v1, v2}.
In this paper, we push the study on projective embeddings of blowup schemes a step
forward. We will investigate the arithmetically Cohen-Macaulayness, the regularity and syzygies of $\vte$. We further examine how these properties behave asymptotically (i.e. when
$t \gg e > 0$).

It was shown by \cite{ch} (see also \cite{v1}) that, when $R$ is Cohen-Macaulay, under certain conditions,
there exists an $f$ such that for all $e > 0$ and $t > fe$, the embedding $\vte$ of
$\tx$ is arithmetically Cohen-Macaulay, i.e. it has a Cohen-Macaulay
coordinate ring. However, no explicit bound for $f$ could be found unless one is in
special situations, for example if it is known that the Rees algebra $\R(I)$ of $I$
is Cohen-Macaulay (see \cite{v1}), or when $\mbox{char } k = 0$, $X = \pp^n$, $Z$
is a scheme of fat points in $X$ and $e=1$ (see \cite[Theorem 2.4]{ge-gi-p}). In the
first section of this paper, we address this problem  again, and generalize the
method of \cite{ge-gi-p} to give an explicit bound for $f$ for any $X$ and $Z$,
when the characteristic of $k$ is 0. More precisely, we give a constant $\delta = \delta(I) > 0$
such that for all $e > 0$ and $t > ed(I) + \delta$ (where $d(I)$ is the maximum
degree of a minimal system of homogeneous generators of $I$), the projective
embedding $\vte$ of $\tx$ is arithmetically Cohen-Macaulay (Theorem \ref{CM}). An
explicit bound for $f$ then can be taken to be $d(I) + \delta(I)$. We also replace the condition on the Cohen-Macaulayness of $R$ (which has been the case in most of previous works) by a slightly weaker condition that $H^i(X, \calo_X) = 0$ for $i = 1, \ldots, n-1$. This turns out to be a necessary condition for the existence of any arithmetically Cohen-Macaulay $\vte$ (Remark \ref{necessary}). Our result is stated as follows.

\begin{thm}[Theorem \ref{CM}]
Suppose $H^i(X, \calo_X) = 0$ for $i = 1, \ldots, n-1$. Then, for $t > ed(I)+\delta$, the projective embedding $\vte$ of $\tx$ is arithmetically
Cohen-Macaulay.
\end{thm}

In the second section of this paper, we study the regularity of the ideal sheaf and
syzygies of $\vte$. Let $\sheaf{J}_{e,t}$ denote the ideal sheaf of $\vte$. We prove that, under a mild condition, the regularity of $\sheaf{J}_{e,t}$
is always bounded above by the dimension of $R$, when $t > ed(I)+\delta$ (Theorem
\ref{reg}). This gives a uniform bound on the shifts of the minimal free resolution
of $\vte$. A very interesting asymptotic behaviour is that, with the additional condition $H^i(X, \calo_X) = 0$ for $i = 1, \ldots, n-1$, the
regularity of $\sheaf{J}_{e,t}$ stabilizes as $t$ gets large compared to
$e$. We also prove this fact in Theorem \ref{reg}.

\begin{thm}[Theorem \ref{reg}]
Suppose $H^0(\tx, w_{\tx}) = 0$, where $w_{\tx}$ is the dualizing sheaf on $\tx$. Then,
\[ \reg \sheaf{J}_{e,t} \le n+1, \ \forall \ t > ed(I)+\delta. \]
If, in addition, $H^i(X, \calo_X) = 0$ for $i = 1, \ldots, n-1$, then for all $t \gg e$,
\[ \reg \sheaf{J}_{e,t} = n+1. \]
\end{thm}

We further restrict our attention to the case when $X$ is a surface and $Z$ is a
0-dimensional subscheme of $X$, and show that for any $p \in \nn$, the embedding
$\vte$ possesses property $N_p$ for all $t \gg e > 0$ (Theorem \ref{np}). This is a
property, introduced by Green and Lazarsfeld (\cite{gr, gl}), that encodes a
lot of information about the scheme.

Results in this paper exhibit and support the philosophy that when embedding a
variety into large projective spaces, it seems to have nicer algebraic properties,
and its algebraic invariants seem to stabilize as the dimension of the ambient space
gets large. This philosophy is the guideline of many recent studies since the
impact of works of Mumford \cite{mum}, its amplification by St. Donat \cite{std},
and the strong theorems and conjectures of Green \cite{gr} and Green and Lazarsfeld
\cite{gl}.

Throughout this paper, our base field $k$ will be algebraically closed and of
characteristic 0.

\section{Cohen-Macaulayness}

Let us first briefly recall some notations and terminology. Unexplained terminology
follow from that of \cite{e} and \cite{h}. Suppose $Y = \proj T \subseteq \pp^l$ is a projective scheme. $Y$ is said to be {\it arithmetically Cohen-Macaulay} if $T$ is a Cohen-Macaulay ring. Suppose $\sheaf{J}$ is the ideal sheaf of $Y$ in $\pp^l$. The {\it regularity} of $\sheaf{J}$, denoted by $\reg \sheaf{J}$, is defined to be the smallest integer $r'$ such that $H^i(\pp^l, \sheaf{J}(r'-i)) = 0$ for all $i > 0$. The regularity can also be interpreted in terms of the minimal free resolution as follows. Let
\[ 0 \rightarrow \bigoplus_j \calo_{\pp^l}(-b_{sj}) \rightarrow \ldots \rightarrow \bigoplus_j \calo_{\pp^l}(-b_{1j}) \rightarrow \bigoplus_j \calo_{\pp^l}(-b_{0j}) \rightarrow \sheaf{J} \rightarrow 0 \]
be the minimal free resolution of $\sheaf{J}$, then $\reg \sheaf{J} = \max\{ b_{ij} - i \ | \ b_{ij} \not= 0 \}$. 

Back to our setup, let $R$ be a finitely generated standard graded $k$-algebra of dimension (n+1), $X = \mbox{Proj } R$ a nonsingular projective scheme. Let $Z$ be a nonsingular closed subscheme of $X$ defined by a homogeneous ideal $I \subseteq R$ of positive height ($I$ saturated). Let $\ti$ be the associated sheaf of $I$ on $X$, and let
$ \pi: \tx = \mbox{Proj } (\oplus_{n \ge 0} \ti^n) \rightarrow X$
be the blowup of $X$ centered at $Z$. Let $E$ be the exceptional divisor of this
blowing up, and $E_0$ the pull-back to $\tx$ of a general hyperplane in $X$. For
each pair of positive integers $e$ and $t$, consider the divisor $\dte = tE_0 - eE$
on $\tx$. Suppose $I$ is generated in degrees at most $d(I)$. By \cite[Lemma
1.1]{ch}, for each $e > 0$ and $t\ge ed(I)+1$, the divisor $\dte$ is very ample on
$\tx$, so from now on, we only consider values of $e > 0$ and $t$ such that $t \ge
ed(I)+1$. For such values of $e$ and $t$, let $\tx \isomap \vte \subseteq
\pp^{N_{e,t}}$ be the projective embedding of $\tx$ given by $\dte$, where $N_{e,t}
= \dim_k H^0(\tx, \dte)-1$. We also denote by $r$ the height of $I$.

Now, suppose $R$ is generated as a $k$-algebra by $x_1, \ldots, x_g \in R_1$, and $I$ is generated by $f_1, \ldots, f_m \in R$, where the degree of $f_j$
is $d_j$ for all $j = 1, \ldots, m$. Let $\R(I) = \oplus_{n \ge 0} I^nt^n$ be the
Rees algebra of $I$. Let $S = k[X_1, \ldots, X_g, Y_1, \ldots, Y_m]$ be a polynomial ring, then $S$ has a natural bi-gradation (induced by $\deg X_i = (1,0), \ \forall i = 1, \ldots, g$, and $\deg Y_j = (d_j, 1), \ \forall j = 1, \ldots, m$). Consider the $k$-algebra homomorphism
\[ S \rightarrow \R(I), \]
given by sending $X_i$ to $x_i$ for all $i = 1, \ldots, g$, and $Y_j$ to $f_jt$ for all $j = 1, \ldots, m$. It is clear that this homomorphism is surjective and gives $\R(I)$ the structure of a bi-graded module over $S$. Suppose
\[ 0 \rightarrow \ldots \rightarrow \oplus_{j}S(-a_{ij}, -b_{ij}) \rightarrow
 \ldots \rightarrow \oplus_{j}S(-a_{0j}, -b_{0j}) \rightarrow \R(I) \rightarrow
0 \] is the bi-graded minimal free resolution of $\R(I)$ over $S$. This minimal
free resolution exists and is finite due to Hilbert's syzygy theorem. Let
\[ c_i = \max_j \{a_{ij} - b_{ij}d(I)\}. \]
Clearly, $\max_i \{ c_i - i \}$ is an invariant of $I$. We denote this invariant by
$c(I)$. Define 
\[ \delta = \delta(I) = \max \{ c(I), 1 \}. \] 

A straight forward adaptation of \cite[Theorem 2.4 and Proposition 4.1]{cht} gives us the following result.

\begin{lemma} \label{cht}
With $\delta = \delta(I)$ defined as above, for any natural number $l$, we have
\[ \reg \sheaf{I}^l \le \reg I^l \le ld(I) + \delta. \]
\end{lemma}

The following proposition is an extension of \cite[Proposition 2.2]{ge-gi-p}.

\begin{pro} \label{projectively_normal}
For $t > ed(I) + \delta$, the projective embedding $\vte$ of $\tx$ is projectively
normal.
\end{pro}

\begin{proof} Suppose $t > ed(I)+\delta$. By \cite[Lemma 1.1]{ch}, $\dte$
is very ample on $\tx$. Let $S_{e,t} = \mbox{Sym}^{*}(H^0(\tx, \dte))$ be the
coordinate ring of $\pp^N = \pp^{N_{e,t}}$ into which $\tx$ is embedded by $\dte$. Let
$I_{e,t} \subseteq S_{e,t}$ be the defining ideal of $\vte$, and $\sheaf{J}_{e,t}$ its associated ideal sheaf in $\pp^N$. From now on, when working with $\vte$, we shall always use these notations of $S_{e,t}, I_{e,t}, \sheaf{J}_{e,t}$ and $N = N_{e,t}$. One has the following exact sequence
\[ 0 \rightarrow I_{e,t} \rightarrow S_{e,t} \rightarrow \oplus_{h \ge 0}
H^0(\tx, \calo_{\tx}(h\dte)) \rightarrow \oplus_{h \ge 0} H^1(\pp^N,
\sheaf{J}_{e,t}(h)) \rightarrow 0. \] 

To prove $\vte$ is projectively normal, it is enough to show that for each
 $h \ge 0$, there is an isomorphism of $k$-vector spaces
\[ \Big(\sfrac{S_{e,t}}{I_{e,t}}\Big)_h \simeq H^0(\tx, \calo_{\tx}(h\dte)), \]
where $\Big(\sfrac{S_{e,t}}{I_{e,t}}\Big)_h$ is the $k$-vector space of monomials of degree $h$ in $\sfrac{S_{e,t}}{I_{e,t}}$.
This is clear for $h = 0$. Suppose $h \ge 1$. It can be observed that a hyperplane in $\pp^N$ when restricted to $\vte$ pulls back to a hypersurface of
degree $t$  in $X$ containing $Z$ with multiplicity at least $e$. Thus,
\begin{eqnarray}
\Big(\sfrac{S_{e,t}}{I_{e,t}}\Big)_h \simeq [(I^e_t)^h]_{ht} = [(I^e)^h]_{ht}  =
(I^{he})_{ht}, \label{eq0}
\end{eqnarray}
where $[(I^e_t)^h]_{ht}$ is the $k$-vector space of monomials of degree $ht$ in the $h$-th power of the ideal generated by $I^e_t$.

Let $\cm = \ti.\calo_{\tx} = \calo_{\tx}(-E)$, then
\[ \calo_{\tx}(h\dte) = \cm^{he} \otimes \calo_{\tx}(htE_0) = \cm^{he} \otimes
 \pi^{*} \calo_X(ht). \]
By \cite[Proposition 10.2]{mat}, we have $\pi_{*} \cm^{he} = \ti^{he}$. Hence,
by the projection formula, we get
\[ \pi_{*} \calo_{\tx}(h\dte) = \pi_{*}(\cm^{he}) \otimes \calo_X(ht) =
 \ti^{he}(ht). \]
Thus,
\begin{eqnarray}
H^0(\tx, \calo_{\tx}(h\dte)) = H^0(X, \pi_* \calo_{\tx}(h\dte)) = H^0(X,
\ti^{he}(ht)) = \dim_k I^{(he)}_{ht}, \label{eq00}
\end{eqnarray}
where $I^{(he)}$ is the saturation ideal of $I^{he}$.

It follows from Lemma \ref{cht} that
\[ \mbox{reg } I^l \le ld(I) + \delta \mbox{ for any } l \ge 1. \]
Thus, for $t > ed(I)+\delta$, we have
\[ ht > hed(I) + h\delta \ge hed(I) + \delta \ge \mbox{reg } I^{he} \ge \mbox{sat } (I^{he}), \]
where $\mbox{sat } (I^{he})$ is the saturation degree of $I^{he}$ (i.e. the least
degree starting from which $I^{he}$ and its saturation ideal $I^{(he)}$ agree).
Therefore, for $t  > ed(I)+\delta$, we have
\[ \dim_k [I^{he}]_{ht} = \dim_k I^{(he)}_{ht}. \]
This, together with (\ref{eq0}) and (\ref{eq00}), proves the proposition.
\end{proof}

The following lemma plays an important role in the study of arithmetically
Cohen-Macaulayness of $\vte$.

\begin{lemma} \label{direct_image}
For $e > 0, t \ge ed(I)+1$, we have
\begin{enumerate}
\item $R^j\pi_{*} \calo_{\tx}(h\dte) = 0$ for all $j > 0$ and $h \ge 0.$
\item $H^i(\tx, \calo_{\tx}(h\dte)) = H^i(X, \ti^{he}(ht))$ for all $i \ge 0$
 and $h \ge 0$.
\end{enumerate}
Here, $\calo_{\tx}$ is the structure sheaf of $\tx$.
\end{lemma}

\begin{proof} Again, we let $\cm = \ti\calo_{\tx} = \calo_{\tx}(-E)$.

{\bf (1)} For each $l \ge 0$, we have
\[ \calo_{\tx}(h\dte) \otimes \pi^{*}\calo_X(l) = \calo_{\tx}(h(tE_0 - eE))
 \otimes \pi^{*}\calo_X(l) = \cm^{he} \otimes \pi^{*}\calo_X(ht+l). \]

Now we will use the argument of \cite[Proposition 1.5]{ch} to get
\[ H^i(\tx, \calo_{\tx}(h\dte) \otimes \pi^{*}\calo_X(l)) = 0, \forall i > 0, l
 \gg 0. \]
Indeed, let $w_{\tx}$ and $w_X$ be the dualizing sheaves on $\tx$ and $X$, respectively. Since the height of $I$ is $r$,
$w_{\tx} = \pi^{*} w_X \otimes \cm^{1-r}$. Thus,
\[ \cm^{he} \otimes \pi^{*}\calo_X(ht+l) = \cm^{he+r-1} \otimes \pi^{*}\calo_X(ht+l) \otimes
(\pi^{*}w_X)^{-1} \otimes w_{\tx}. \] Observe that for any $c,d \in \nn$, $\cm^{d}
\otimes \pi^{*}\calo_X(c) = I^d_c.\calo_{\tx}$. It is known that $w_X^{-1}(g)$ is
very ample on $X$ for some $g > 0$. Thus, by \cite[Lemma 1.1]{ch}, for $l >
(he+r-1)d(I)+g - ht$ (i.e. $ht+l > (he+r-1)d(I) + g$), the divisor $\cm^{he+r-1}
\otimes \pi^{*}\calo_X(ht+l) \otimes (\pi^{*}w_X)^{-1}$ is very ample on $\tx$. By Kodaira's vanishing theorem, one has for $l > (he+r-1)d(I) +  g - ht$ and $i > 0$,
\begin{eqnarray}
H^i(\tx, \calo_{\tx}(h\dte) \otimes \pi^{*}\calo_X(l)) & = & 0. \label{eq1}
\end{eqnarray}

Next, for each $i > 0$ and $h \ge 0$, since $\calo_X(1)$ is very ample on $X$, by
Serre's theorem, there exists an integer $l_{i,h}$ such that for all $l \ge
l_{i,h}, p > 0$ and $q \le i$,
\[ H^p(X, R^q\pi_{*} \calo_{\tx}(h\dte) \otimes \calo_X(l)) = 0. \]
By the projection formula, we then have, for all $l \ge l_{i,h}, p > 0$ and $q
\le i$,
\begin{eqnarray}
H^p(X, R^q\pi_{*}(\calo_{\tx}(h\dte) \otimes \pi^{*} \calo_{X}(l))) & = & 0.
 \label{eq2}
\end{eqnarray}

Consider the Leray spectral sequence
\[ E^{p,q}_2 = H^p(X, R^q\pi_{*}(\calo_{\tx}(h\dte) \otimes \pi^{*}
 \calo_{X}(l))) \Rightarrow
H^{p+q}(\tx, \calo_{\tx}(h\dte) \otimes \pi^{*} \calo_{X}(l)). \]

By (\ref{eq2}), this Leray spectral sequence (when used to compute $H^i(\tx,
\calo_{\tx}(h\dte) \otimes \pi^{*}\calo_X(l))$) concentrates on its vertical
 boundary. This and (\ref{eq1}) imply that, for all $l > \max\{l_{i,h},
(he+r-1)d(I) + g\}$ and all  $i > 0$,
\begin{eqnarray}
\Gamma(X, R^i\pi_{*}(\calo_{\tx}(h\dte) \otimes \pi^{*} \calo_{X}(l))) & = &
 H^i(\tx,
\calo_{\tx}(h\dte) \otimes \pi^{*} \calo_{X}(l)) = 0. \label{eq3}
\end{eqnarray}
Since $\calo_X(1)$ is very ample on $X$, we also know that for $l \gg 0$,
\[ R^i\pi_{*}(\calo_{\tx}(h\dte) \otimes \pi^{*}
\calo_{X}(l)) = R^i\pi_{*} \calo_{\tx}(h\dte) \otimes \calo_X(l) \]
is generated by global sections. Thus, (\ref{eq3}) implies
$R^i\pi_{*}(\calo_{\tx}(h\dte)) = 0, \ \forall i > 0$. The result is proved.

{\bf (2)} It follows from the Leray spectral sequence
\[ E^{p,q}_2 = H^p(X, R^q\pi_{*} \calo_{\tx}(h\dte)) \Rightarrow H^{p+q}(\tx,
\calo_{\tx}(h\dte)), \]
and the result in part (1) that
\[ H^i(X, \pi_{*} \calo_{\tx}(h\dte)) = H^i(\tx, \calo_{\tx}(h\dte)), \ \forall i
 \ge 0 \mbox{ and } h \ge 0. \]
We also have, by the projection formula and \cite[Proposition 10.2]{mat},
\[ \pi_{*} \calo_{\tx}(h\dte) = \pi_{*} (\cm^{he} \otimes \pi^{*}\calo_X(ht)) =
\pi_{*} \cm^{he} \otimes \calo_X(ht) = \ti^{he}(ht). \]
Thus, the result follows, and the lemma is proved.
\end{proof}

The main result of this section is stated as follows.

\begin{thm} \label{CM}
Suppose $H^i(X, \calo_X) = 0$ for $i = 1, \ldots, n-1$. Then, for $t > ed(I)+\delta$, the projective embedding $\vte$ of $\tx$ is arithmetically
Cohen-Macaulay.
\end{thm}

\begin{proof}
Having Lemma \ref{direct_image}, the proof now follows in the same line as that
 of \cite[Theorem 2.4]{ge-gi-p} with a slight modification.

Suppose $t > ed(I)+\delta$. As in Proposition \ref{projectively_normal}, we
 let $S_{e,t} = \mbox{Sym}^{*}H^0(\tx, \dte)$ be the coordinate ring of $\pp^N =
\pp^{N_{e,t}}$, let $I_{e,t} \subseteq S_{e,t}$ be the defining ideal of
 $\vte$, and $\sheaf{J}_{e,t}$ its ideal sheaf in $\pp^N$. Since $\dim \vte = \dim X = n$, to prove $\vte$ is arithmetically Cohen-Macaulay, it is enough to show that $H^1(\pp^N,
\sheaf{J}_{e,t}(h)) = 0$  for all $h \in \zz$ and $H^i(\pp^N, \calo_{\vte}(h)) =
0$ for all $i = 1, \ldots, n-1$ and for all $h \in \zz$.

By Proposition \ref{projectively_normal}, we know $H^1(\pp^N,
\sheaf{J}_{e,t}(h)) = 0$ for all $h \ge 0$. For $h < 0$, it is clear that
$\calo_{\vte}(-h)$ is a very ample invertible sheaf on $\vte$. Therefore, by
Kodaira's vanishing theorem, we get $H^0(\pp^N, \calo_{\vte}(h)) = 0$ for all $h <
0$. This implies that $H^1(\pp^N, \sheaf{J}_{e,t}(h)) = 0$ for all $h < 0$. Thus,
\[ H^1(\pp^N, \sheaf{J}_{e,t}(h)) = 0, \mbox{ for all } h \in \zz. \]

Let us now consider $H^i(\pp^N, \calo_{\vte}(h))$ ($i = 1, \ldots, n-1$) for $h \ge 0$. By Lemma \ref{direct_image}, we have, for all $i = 1, \dots, n-1$,
\[ H^i(\pp^N, \calo_{\vte}(h)) = H^i(\tx, \calo_{\tx}(h\dte)) = H^i(X, \ti^{he}(ht)). \]
Moreover, it follows from Lemma \ref{cht} that, for $h > 0$,
\[ \mbox{reg } \ti^{he} \le hed(I) + \delta \le hed(I)+h\delta < ht. \]
Thus, $H^i(X, \ti^{he}(ht)) = 0$ for all $i = 1, \dots, n-1$, and $h > 0$. Together with the hypothesis, we have
\[ H^i(\pp^N, \calo_{\vte}(h)) = 0, \mbox{ for } i = 1,\ldots, n-1, \mbox{ and } h \ge 0. \]

For $h < 0$, we can again use the Kodaira's vanishing theorem to obtain
\[ H^i(\pp^N, \calo_{\vte}(h)) = 0, \mbox{ for } i = 1, \ldots, n-1. \]
This is due to the fact that $\calo_{\vte}(-h)$ is a very ample invertible sheaf on $\vte$ for all $h < 0$. The theorem is proved.
\end{proof}

\begin{cor}
Suppose $R$ is Cohen-Macaulay. Suppose $I$ is generated by $s$ elements and the Rees algebra $\R(I)$ of $I$ is Cohen-Macaulay. Then, for all $e > 0$ and $t > ed(I) + (s-1)(d-1)$, the projective
embedding $\vte$ of $\tx$ is arithmetically Cohen-Macaulay.
\end{cor}

\begin{proof}
It follows from \cite[Proposition 4.1]{v1} that, in the minimal free resolution
 of the Rees algebra $\R(I)$ of $I$, the bi-graded Betti numbers satisfy
\[ a_{ij} \le sd - (s-1)+i. \]
Thus, $\delta(I) \le (s-1)[d(I)-1]$. The result now follows from Theorem
 \ref{CM}.
\end{proof}

\begin{remark} \label{necessary} 
It follows from \cite[Theorem 1.3]{tc} that the condition $H^i(X, \calo_X) = 0$ for all $i = 1, \ldots, n-1$ is necessary. This is because an arithmetically Cohen-Macaulay $\vte$ would yield an {\it arithmetic Macaulayfication} of $X$ (see \cite{tc} for the definition and further results on arithmetic Macaulayfication of projective schemes).
\end{remark}

\begin{remark} We should point out that with exactly the same proofs, all our results
in this section are true even if $Z$ is only a locally complete intersection
subscheme of $X$.
\end{remark}

\section{Regularity and $N_p$ property}

In this section, we study how the regularity and the syzygies of $\vte$ behave
asymptotically. We use the same notations and terminology as in the previous
section.

Let us also recall the notions of Koszul complex and Koszul cohomology which
 were studied in \cite{gr}.

Let $Y$ be a projective scheme. Let $\cl$ be a very ample line bundle and $\cf$
 a coherent sheaf on $Y$. Let $W = H^0(Y, \cl)$ and $S = \mbox{Sym}^{*} W$.
Then,  $S$ is the homogeneous coordinate ring of $\pp(W)$, the projective space into
which $\cl$ embeds $Y$. Let $B = B(\cf, \cl) = \oplus_{q \in \zz} H^0(Y, \cf \otimes
q\cl) = \oplus_{q \in \zz} B_q$ a $S$-graded module.

\begin{mydef} The {\it Koszul complex} of $B$ is defined to be
\[ \ldots \rightarrow \wedge^{p+1}W \otimes B_{q-1}
 \stackrel{d_{p+1,q-1}}{\rightarrow} \wedge^pW \otimes B_q
\stackrel{d_{p,q}}{\rightarrow} \wedge^{p-1}W \otimes  B_{q+1} \rightarrow
\ldots \] and the {\it Koszul cohomology groups} of $B$ are defined to be
\[ \ck_{p,q}(\cf, \cl) = \sfrac{\mbox{ker } d_{p,q}}{\mbox{im } d_{p+1,q-1}},
 \mbox{ for } p,q \in \zz. \]
\end{mydef}

When $\cl = \cl(D)$ is the invertible sheaf corresponding to a divisor $D$ and
 $\cf = \tilde{M}$ is the sheaf associated to a module $M$, we write
$\ck_{p,q}(M, D)$  for $\ck_{p,q}(\cf, \cl)$. When $\cf$ is the structure sheaf
$\calo_Y$, we also write $\ck_{p,q}(\cl)$ for $\ck_{p,q}(\cf, \cl)$.

Recall further that $\sheaf{J}_{e,t}$ denotes the ideal sheaf of $\vte$. The following theorem gives an upper bound for the regularity of $\sheaf{J}_{e,t}$, and shows that this regularity stabilizes for $t \gg e$. 

\begin{thm} \label{reg}
Suppose $H^0(\tx, w_{\tx}) = 0$, where $w_{\tx}$ is the dualizing sheaf on $\tx$. Then,
\[ \reg \sheaf{J}_{e,t} \le n+1, \ \forall \ t > ed(I)+\delta. \]
If, in addition, $H^i(X, \calo_X) = 0$ for $i = 1, \ldots, n-1$, then for all $t \gg e$,
\[ \reg \sheaf{J}_{e,t} = n+1. \]
\end{thm}

\begin{proof} Let $S_{e,t}$ and $I_{e,t}$ be as before (see Proposition
\ref{projectively_normal}). It follows from Proposition \ref{projectively_normal} that for $t > ed(I)+\delta$,
\[ \sfrac{S_{e,t}}{I_{e,t}} \simeq \oplus_{h \ge 0} H^0(\tx,
\calo_{\tx}(h\dte)). \] 
Thus, by Green's syzygy theorem (\cite[1.b.4]{gr}), the minimal free resolution of $S_{e,t}/I_{e,t}$ is
given by
\[ 0 \rightarrow F_s \rightarrow \ldots \rightarrow F_1 \rightarrow
F_0 =  S_{e,t} \rightarrow S_{e,t}/I_{e,t} \rightarrow 0, \]
for some $s \ge N-n$ and 
\[ F_i = \bigoplus_{q \ge 1} \ck_{i,q}(\dte) \otimes S(-i-q), \mbox{ for }
i =1,  \ldots, s. \]

We first prove $\mbox{reg } \sheaf{J}_{e,t} \le n+1$ for all $t > ed(I)+\delta$.
This is equivalent to showing that $\ck_{i, q}(\dte) = 0$ for all $1 \le i \le s$
and $q \ge n+1$, when $t > ed(I)+\delta$. Consider $1 \le i \le s$ and $q \ge
n+1$. Let $K_{\tx}$ be the canonical divisor of $\tx$, then $w_{\tx}$ is the sheaf
associated to $K_{\tx}$ on $\tx$. It follows from the proof of Theorem \ref{CM} that
\[ H^i(\tx, \calo_{\tx}(h\dte)) = H^i(X, \ti^{he}(ht)) = 0, \]
for all $i=1, \ldots, n-1$, and $h > 0$. Thus, Green's Duality theorem (\cite[2.c.6]{gr}) gives us
\[ \ck_{i,q}(\dte)^{*} \simeq \ck_{N-n-i, n+1-q}(K_{\tx}, \dte), \ \forall \ q \ge n+1. \]
This implies that, if $i > N-n$, then $\ck_{i,q}(\dte) = 0$. Suppose $i \le N-n$. By Green's Vanishing theorem (\cite[3.a.1]{gr}), it suffices now to show that
\begin{eqnarray}
h^0(\tx, w_{\tx} \otimes \calo_{\tx}((n+1-q)\dte)) & \le & N-n-i, \ \forall \ t >
ed(I)+\delta. \label{eq13}
\end{eqnarray}
This is indeed true. If $q = n+1$, this follows from the hypothesis that $H^0(\tx,
w_{\tx}) = 0$. Otherwise, if $q > n+1$, then by Serre's duality one has $H^0(\tx,
w_{\tx} \otimes \calo_{\tx}((n+1-q)\dte)) = H^n(\tx, \calo_{\tx}((q-n-1)\dte))$.
Thus, by Lemma \ref{direct_image} and \cite[Proposition 10.2]{mat}, we have
\[ H^0(\tx, w_{\tx} \otimes \calo_{\tx}((n+1-q)\dte)) = H^n(X, \ti^{e(q-n-1)}(t(q-n-1))). \]
It follows from Lemma \ref{cht} that for $t > ed(I)+\delta$,
\[ \mbox{reg } \ti^{e(q-n-1)} < t(q-n-1), \]
whence
\[ H^n(X, \ti^{e(q-n-1)}(t(q-n-1))) = 0. \]
Therefore,
\[ H^0(\tx, w_{\tx} \otimes \calo_{\tx}((n+1-q)\dte)) = 0, \ \forall \ t > ed(I)+\delta, \]
and (\ref{eq13}) is proved.

Suppose, the additional condition $H^i(X, \calo_X) = 0$ for $i = 1, \ldots, n-1$ is satisfied. This, together with the proof of Theorem \ref{CM}, shows that 
\[ H^i(\tx, \calo_{\tx}(h\dte)) = H^i(X, \ti^{he}(ht)) = 0, \]
for all $i=1, \ldots, n-1$, and $h \ge 0$. Thus, Green's Duality theorem (\cite[2.c.6]{gr}) holds for $\ck_{N-n, n}(\dte)$, and we have
\begin{eqnarray}
\ck_{N-n,n}(\dte)^{*} & \simeq & \ck_{0, 1}(K_{\tx}, \dte). \label{eq8}
\end{eqnarray}
Let $W = H^0(\tx, \calo_{\tx}(\dte))$, the Koszul complex of $\oplus_{h \ge 0}
H^0(\tx, \calo_{\tx}(h\dte))$ at degree $(0,1)$ is
\[ \wedge W \otimes H^0(\tx, w_{\tx}) \stackrel{d_{1,0}}{\rightarrow} W \otimes
H^0(w_{\tx} \otimes \calo_{\tx}(\dte)) \stackrel{d_{0,1}}{\rightarrow} 0. \] By
 the hypothesis, since $H^0(\tx, w_{\tx}) = 0$, we have
\begin{eqnarray}
\ck_{0,1}(K_{\tx}, \dte) & = & W \otimes H^0(w_{\tx} \otimes \calo_{\tx}(\dte)).
 \label{eq9}
\end{eqnarray}

Recall that $r$ is the height of $I$, $w_X$ is the dualizing sheaf on $X$, and $\cm
= \ti \calo_{\tx} = \calo_{\tx}(-E)$, we have
\[ H^0(\tx, w_{\tx} \otimes \calo_{\tx}(\dte)) = H^0(X, \pi_*(w_{\tx} \otimes
\calo_{\tx}(h\dte))) = H^0(\tx, \pi_{*}(\cm^{e+1-r} \otimes \pi^{*} w_X(t))). \]
By the projection formula and by \cite[Proposition 10.2]{mat}, we now obtain
\[ H^0(\tx, w_{\tx} \otimes \calo_{\tx}(\dte)) = H^0(X, \ti^{e+1-r} \otimes w_X \otimes
\calo_X(t)), \]
where it is understood that $\ti^{e+1-r} = \calo_X$ if $e \le r-1$.

It can be seen that $\calo_X(1)$ is an ample invertible sheaf on $X$, so for $t
\gg e$, $\ti^{e+1-r} \otimes w_X \otimes \calo_X(t)$ is generated by its global
sections. Therefore, since $\ti^{e+1-r} \otimes w_X \otimes \calo_X(t)$ cannot be
the zero sheaf on $X$, we must have $H^0(X, \ti^{e+1-r} \otimes w_X \otimes
\calo_X(t)) \not= 0$ for all $t \gg e$. It now follows from (\ref{eq8}) and
(\ref{eq9}) that, for all $t \gg e$,
\[ \ck_{0,1}(K_{\tx}, \dte) \not= 0, \mbox{ i.e. } \ck_{N-n, n}(\dte) \not= 0. \]
This implies that $\reg \sheaf{J}_{e,t} \ge n+1$. The result now follows. The
theorem is proved.
\end{proof}

\begin{remark} \label{rem}
One can observe that even without the condition $H^0(\tx, w_{\tx}) = 0$, the proof of Theorem \ref{reg} still shows $\mbox{reg } \sheaf{J}_{e,t} \le n+1$ for all $t \gg e > 0$. This is because for $t \gg e$, $N = N_{e,t} \gg 0$.
\end{remark}

Before moving on to study the syzygies of $\vte$, let us recall the notion of {\it
property $N_p$} (from \cite{gl}).

\begin{mydef} Let $Y$ be a projective variety and let $\cl$ be a very ample line bundle on $Y$
defining an embedding $\varphi_{\cl}: Y \hookrightarrow \pp = \pp(H^0(Y, \cl)^{*})$.
\begin{enumerate}
\item The line bundle $\cl$, or the embedding $\varphi_{\cl}(Y)$ of $Y$, is said to have
{\it property} $N_0$ if $\varphi_{\cl}(Y)$ is projectively normal, i.e. $\cl$ is
normally generated.
\item Let $S = \mbox{Sym}^{*}H^0(Y, \cl)$, the homogeneous coordinate ring of the projective
space $\pp$. Suppose $A$ is the homogeneous coordinate ring of $\varphi_{\cl}(Y)$ in
$\pp$, and
\[ 0 \rightarrow F_n \rightarrow F_{n-1} \rightarrow \ldots \rightarrow F_0 \rightarrow A
\rightarrow 0 \]
is a minimal free resolution of $A$. The line bundle $\cl$, or the embedding
$\varphi_{\cl}(Y)$ of $Y$, is said to have {\it property} $N_p$ (for $p \in \nn$) if and
only if it has property $N_0$, $F_0 = S$ and $F_i = S(-i-1)^{\alpha_i}$ with $\alpha_i
\in \nn$ for all $1 \le i \le p$.
\end{enumerate}
\end{mydef}

In what comes next, we restrict our attention to the case when $X$ is a surface,
i.e. $\dim R = 3$, and $Z$ is a nonsingular 0-dimensional subscheme of $X$. We will
show that for any $p \in \nn$, the projective embedding $\vte$ of $\tx$ possesses
property $N_p$ for all $t \gg e > 0$. We apply the method which was demonstrated in
\cite{n2}.

\begin{thm} \label{np}
Suppose $X$ is a nonsingular surface such that $H^1(X, \calo_X) = 0$, and $Z$ is a nonsingular 0-dimensional
subscheme of $X$. Let $\tx$ be the blowup of $X$ centered at $Z$. Then, for any $p
\in \nn$, the projective embedding $\vte$ of $\tx$ possesses property $N_p$ for all $t
\gg e > 0$.
\end{thm}

\begin{proof} For $p = 0$, the theorem is proved in Proposition \ref{projectively_normal}.
Suppose now that $p > 0$. Let $S = S_{e,t}, I_{e,t}$ and $\sheaf{J}_{e,t}$ be as
before. It follows from Theorem \ref{reg} and Remark \ref{rem} that for all $t \gg e$,
\[ \mbox{reg } \sheaf{J}_{e,t} \le 3. \]

It also follows from Theorem \ref{CM} that for all $t \gg e$, $\vte$ is
arithmetically Cohen-Macaulay. Thus, for all $t \gg e$, the defining ideal $I_{e,t}$ of $\vte$ has the following minimal free resolution
\[ \begin{array}{ccccccccccccc}
& & S(-N)^{\beta_{N-2,N}} & & & & S(-4)^{\beta_{2,4}} & & S(-3)^{\beta_{1,3}}
& & & & \\
& & \bigoplus & & & & \bigoplus & & \bigoplus & & & & \\
0 & \rightarrow & S(-(N-1))^{\beta_{N-2, N-1}} & \rightarrow & \ldots & \rightarrow &
S(-3)^{\beta_{2,3}} & \rightarrow & S(-2)^{\beta_{1,2}} & \rightarrow & I_{e,t} &
\rightarrow & 0, \end{array} \]
where $N = N_{e,t} = \dim_k H^0(\tx, \dte) -1$. To prove $\vte$ possesses
property $N_p$, it is enough to show that $\beta_{i, i+2} = 0$ for all $1 \le i \le
p$.

Let $\cc$ be a general hyperplane section of $\vte$, then $\cc$ is an arithmetically
Cohen-Macaulay curve in $\pp^{N-1}$ with the same minimal free resolution as that of
$\vte$. Let $A$ be the homogeneous coordinate ring of $\cc$ in $\pp^{N-1}$. Then,
the Betti numbers of $\cc$ (also of $\vte$) are given by
\[ \beta_{i,j} = \mbox{Tor}^T_{i}(A, k)_j \ \forall i, j, \]
where $T = k[x_0, \ldots, x_{N-1}]$ is the coordinate ring of $\pp^{N-1}$.

Now, let $i(E) = E.E$ be the self-intersection number of $E$. Let $A_{e,t} =
S_{e,t}/I_{e,t}$ be the homogeneous coordinate ring of $\vte$, and $\h_{\vte}$ the
Hilbert function of $\vte$, i.e.
\[ \h_{\vte}(\lambda) = \dim_k (A_{e,t})_{\lambda}. \]
For $t \gg e$, since $\vte$ is arithmetically Cohen-Macaulay, we have
\[ \dim_k (A_{e,t})_{\lambda} = \dim_k H^0(\vte, \calo_{\vte}(\lambda)) =
\dim_k H^0(\tx, \calo_{\tx}(\lambda \dte)). \]
Since $\vte$ is arithmetically Cohen-Macaulay, we also have $H^1(\tx,
\calo_{\tx}(\lambda
\dte)) = 0$ for all $\lambda \in \zz$. Thus, by the Riemann-Roch theorem, we have
\[ \dim_k (A_{e,t})_{\lambda} = h^0(\tx, \lambda \dte) \le \sfrac{1}{2}[(\lambda \dte)^2 -
\lambda \dte. K_{\tx}] + 1 + p_{\tx}, \]
where $K_{\tx}$ and $p_{\tx}$ are the canonical divisor and the arithmetic genus of
$X$. Let $\delta_{e,t}$ be the degree of $\vte$, then we now have
\[ \delta_{e,t} = \mbox{deg } \vte \le (\dte)^2 = (tE_0 - eE)^2 = t^2 + e^2i(E). \]

Let $\x = \cc \cap H$ be a general hyperplane section of $\cc$, then $\x$ is a set
of $\delta_{e,t}$ points in $\pp^{N-2}$. Let $\sheaf{J}_C$ be the ideal sheaf of
$\cc$ in $\pp^{N-1}$, $\sheaf{J}_{\x}$ the ideal sheaf of $\x$ in $H \simeq
\pp^{N-2}$, and $\h_{\x}$ the Hilbert function of $\x$. Since $\cc$ is
arithmetically Cohen-Macaulay, it follows from the exact sequence
\[ 0 \rightarrow \sheaf{J}_C \rightarrow \sheaf{J}_C(1) \rightarrow \sheaf{J}_{\x}(1)
\rightarrow 0 \]
that
\[ 0 \rightarrow H^1(\pp^{N-2}, \sheaf{J}_{\x}(1)) \rightarrow H^2(\pp^{N-1}, \sheaf{J}_C)
\rightarrow H^2(\pp^{N-1}, \sheaf{J}_C(1)) \rightarrow 0. \]
This implies
\begin{eqnarray}
h^2(\pp^{N-1}, \sheaf{J}_C) - h^2(\pp^{N-1}, \sheaf{J}_C(1)) & = & h^1(\pp^{N-2},
\sheaf{J}_{\x}(1)) = \delta_{e,t} - \h_{\x}(1). \label{eq5}
\end{eqnarray}

Let $Z_e$ be the subscheme of $X$ defined by $I^e$. Then, $Z_e$ is also a
0-dimensional subscheme of $X$. We, therefore, know that the Hilbert function of
$Z_e$ eventually stabilizes at $\mbox{deg } Z_e$. Since $R$ is a standard graded
$k$-algebra (i.e. $R$ is generated as a $k$-algebra by $R_1$) and $\dim R = 3$, we
have $\dim_k R_t \ge {t+2 \choose 2}$. Thus, for $t \gg e$,
\[ N = \dim_k R_t - \mbox{deg } Z_e \ge {t+2 \choose 2} - \mbox{deg } Z_e. \]
Clearly, for $t \gg e$ (more precisely, when $t \ge \sfrac{e^2i(E)+2\deg
Z_e+1+p}{3}$),
\[ 2N-3 \ge t^2 + e^2i(E) +p \ge \delta_{e,t}+p. \]
By \cite{b}, we now have, for $t \gg e$,
\begin{eqnarray}
\h_{\x}(1) & \ge & \min \{\delta_{e,t}, N-1\} \ge \delta_{e,t} - (N-2)+p. \label{eq6}
\end{eqnarray}
(\ref{eq5}) then implies
\[ h^2(\pp^{N-1}, \sheaf{J}_C) - h^2(\pp^{N-1}, \sheaf{J}_C(1)) \le (N-2) - p < N-2. \]
By \cite[Corollary 3.3]{n}, this only happens if $h^2(\pp^{N-1}, \sheaf{J}_C(1)) =
0$. Thus, for $t \gg e$, (\ref{eq5}) and (\ref{eq6}) give
\[ h^2(\pp^{N-1}, \sheaf{J}_C) = h^1(\pp^{N-2}, \sheaf{J}_{\x}(1)) = \delta_{e,t} -
\h_{\x}(1) \le (N-2) - p. \]

Now, let $w_A$ be the canonical module of $A$. Then,
\[ \dim_k (w_A)_0 = \dim_k [H^2_{\goth{m}}(A)]_0 = h^1(\pp^{N-1}, \calo_{\cc}) =
h^2(\pp^{N-1}, \sheaf{J}_C) \le (N-2)-p. \]

Furthermore, by \cite[Lemma 2.4]{n2}, $w_A$ is torsion free as an $A$-module. Thus,
the vanishing theorem of \cite[Theorem 1.1]{ek} gives us, for $t \gg e$,
\[ [\mbox{Tor}^T_s(w_A, k)]_s = 0, \ \forall \ s \ge (N-2)-p. \]
By duality, we now obtain, for $t \gg e$,
\[ [\mbox{Tor}^T_i(A,k)]_{i+2} = 0, \ \forall \ i \le p, \]
i.e. for $t \gg e$,
\[ \beta_{i, i+2} = 0 \ \forall \ i \le p. \]
This implies that $\vte$ possesses property $N_p$. The theorem is proved.
\end{proof}

{\bf Remark:} It would be interesting to have an explicit bound $d_{e,p}$ such that
for any $p \in \nn$ and $e > 0$, the projective embedding $\vte$ of $\tx$ has
property $N_p$ for all $t \ge d_{e,p}$. It follows from Theorem \ref{CM} and the
proof of Theorem \ref{np} that a bound could be given as
\[ d_{e,p} = \max \Big\{ e[d(I)+\delta(I)]+1, \sfrac{e^2i(E)+2\deg Z_e+1+p}{3} \Big\}, \]
i.e. for any $p \in \nn$ and $e > 0$, $\vte$ has property $N_p$ for all
\[ t \ge \max \Big\{ e[d(I)+\delta(I)]+1, \sfrac{e^2i(E)+2\deg Z_e+1+p}{3} \Big\}. \]

{\it Acknowledgement}. {\sf The author would like to thank Prof. S. D. Cutkosky for
many useful discussions on topics related to the materials of this paper. The author would also like to thank the referees for helpful suggestions making the paper more readable}.

\end{document}